\input amstex
\documentstyle{amsppt}
\pagewidth{5.4in}
\pageheight{7.6in}
\magnification=1200
\TagsOnRight
\NoRunningHeads
\topmatter
\title
\bf Removable singularity of the polyharmonic equation
\endtitle
\author
Shu-Yu Hsu
\endauthor
\affil
Department of Mathematics\\
National Chung Cheng University\\
168 University Road, Min-Hsiung\\
Chia-Yi 621, Taiwan, R.O.C.\\
e-mail:syhsu\@math.ccu.edu.tw
\endaffil
\date
Feb 7, 2007
\enddate
\address
e-mail address:syhsu\@math.ccu.edu.tw
\endaddress
\abstract
Let $x_0\in\Omega\subset\Bbb{R}^n$, $n\ge 2$, be a domain and let $m\ge 2$. 
We will 
prove that a solution $u$ of the polyharmonic equation $\Delta^mu=0$ in 
$\Omega\setminus\{x_0\}$ has a removable singularity at $x_0$ if and only 
if $|\Delta^ku(x)|=o(|x-x_0|^{2-n})\quad\forall k=0,1,2,\dots,m-1$ as 
$|x-x_0|\to 0$ for $n\ge 3$ and $=o(\log (|x-x_0|^{-1}))\quad\forall 
k=0,1,2,\dots,m-1$ as $|x-x_0|\to 0$ for $n=2$. For $m\ge 2$  we will 
also prove that $u$ has a removable singularity at $x_0$ if 
$|u(x)|=o(|x-x_0|^{2m-n})$ as $|x-x_0|\to 0$ for $n\ge 3$ and $|u(x)|
=o(|x-x_0|^{2m-2}\log (|x-x_0|^{-1}))$ as $|x-x_0|\to 0$ for $n=2$.
\endabstract
\endtopmatter
\NoBlackBoxes
\define \pd#1#2{\frac{\partial #1}{\partial #2}}
\define \1{\partial}
\define \2{\overline}
\define \3{\varepsilon}
\define \4{\widetilde}
\define \5{\underline}
\define \6{\hat}
\document

Recently there is a lot of study of the singularities of equations.
Singularities of solutions of the Ricci flow was studied by R.S.~Hamilton
\cite{H} and G.~Perelman \cite{P1},\cite{P2}. Removable singularities of the
solution of the Ricci flow equation on $\Bbb{R}^2$ was studied by S.Y.~Hsu
in \cite{Hs1}. Removable singularities of the solution of the harmonic map 
assoicated with the standard solution of the Ricci flow and removable 
singularities of the heat equation was studied by S.Y.~Hsu in
\cite{Hs2}. 

Let $x_0\in\Omega\subset\Bbb{R}^n$, $n\ge 2$, be a domain.
It is well known \cite{F} that a solution $u$ of the harmonic equation
$$
\Delta u=0\quad\text{ in }\Omega\setminus\{x_0\}
$$
has a removable singularity at $x_0$ if and only if
$$
|u(x)|=\left\{\aligned
&o(|x-x_0|^{2-n})\quad\text{ as }x\to x_0\quad\text{ if }n\ge 3\\
&o(\log |x-x_0|)\quad\text{ as }x\to x_0\quad\text{ if }n=2.
\endaligned\right.
$$
In this paper we will generalize the above result to the case of polyharmonic
equation. We refer the reader to the book \cite{ACL} by N.~Aronszajn,
T.M.~Creese and L.J.~Lipkin for various properties of polyharmonic functions
and the papers \cite{A}, \cite{CGS}, \cite{PV}, \cite{X},  by D.H.~Armitage,
L.A.~Cafferelli, B.~Gidas and J.~Spruck, J.Pipher and G.C.~Verchota, and 
X.~Xu, for some recent related results in this direction.

We now state some definitions. For any $R>0$ and $x_0\in\Bbb{R}^n$, let 
$B_R(x_0)=\{x\in\Bbb{R}^n:|x-x_0|<R\}$ and $B_R=B_R(0)$. Let $G(x,y)$ be 
the Green function for the Laplacian on unit ball $B_1$. That is 
$\Delta_xG(x,y)=-\delta_y(x)$ and $G(x,y)=0$ for all $x\in\1 B_1$, $|y|<1$, 
where $\delta_y$ is the delta mass at $y$. 

We say that a solution of of the polyharmonic equation
$$
\Delta^mu=0\tag 1
$$
in $\Omega\setminus\{x_0\}$ has a removable singularity at $x_0$ if there
exists a constant $\delta>0$ such that $\2{B_{\delta}(x_0)}\subset\Omega$
and a smooth solution $v$ of (1) in $B_{\delta}(x_0)$ such that
$u(x)=v(x)$ for all $0<|x-x_0|<\delta$.  

We choose a radially symmetric function $\phi\in C_0^{\infty}(\Bbb{R}^n)$, 
$0\le\phi\le 1$, $\int_{\Bbb{R}^n}\phi\,dx=1$, such that $\phi (x)=1$ for 
any $|x|\le 1/2$ and $\phi (x)=0$ for any $|x|\ge 1$. For any $\3>0$ and 
function $f$, let $\phi_{\3}(x)=\3^{-n}\phi (x/\3)$ and
$$
\phi_{\3}\ast f(x)=\frac{1}{\3^n}\int_{|z|\le\3}f(x-z)\phi\biggl (
\frac{z}{\3}\biggr )\,dz.
$$
We will assume $m\ge 2$ and $n\ge 2$ for the rest of the paper.

\proclaim{\bf Theorem 1}
Let $x_0\in\Omega$. A solution $u$ of (1) in $\Omega\setminus\{x_0\}$ has a 
removable singularity at $x_0$ if and only if 
$$
|\Delta^ku(x)|=\left\{\aligned
&o(|x-x_0|^{2-n})\qquad\forall k=0,1,2,\dots,m-1\quad\text{ as }
|x-x_0|\to 0\quad\text{ if }n\ge 3\\
&o(\log(1/|x-x_0|))\quad\forall k=0,1,2,\dots,m-1\quad\text{ as }|x-x_0|\to 0
\quad\text{ if }n=2.
\endaligned\right.\tag 2
$$
\endproclaim
\demo{Proof}
Suppose $u$ has a removable singularity at $x_0$. Then there exists 
$\2{B_{\delta}(x_0)}\subset\Omega$ and a smooth solution $v$ of (1) in 
$B_{\delta}(x_0)$ such that $u(x)=v(x)$ for all $0<|x-x_0|<\delta$. Hence 
(2) holds.

Conversely suppose (2) holds. Without loss of generality we may assume that
$x_0=0$ and $\2{B}_1\subset\Omega$. For any $x\in\2{B}_1$, let
$$
v_1(x)=\int_{\1 B_1}\frac{\1 G}{\1 n_y}(x,y)\Delta^{m-1}u(y)\,d\sigma (y)
\tag 3
$$ 
and
$$
v_i(x)=\int_{B_1}G(x,y)v_{i-1}(y)\,dy
+\int_{\1 B_1}\frac{\1 G}{\1 n_y}(x,y)\Delta^{m-i}u(y)\,d\sigma (y)
\quad\forall i=2,3,\dots,m\tag 4
$$
where $\1/\1 n_y$ is the derivative with respect to the unit inward normal 
at the boundary $y\in\1 B_1$ of the domain $B_1$. Then $v_1\in C^{\infty}
(B_1)\cap C(\2{B}_1)$ and
$$\left\{\aligned
&\Delta v_1=0\qquad\qquad\qquad\text{in }B_1\\
&v_1(x)=\Delta^{m-1}u(x)\quad\text{ on }\1 B_1.\endaligned\right.\tag 5
$$
By (4) $v_2\in C^{\infty}(B_1)\cap C(\2{B}_1)$. Then by (4) and an induction
argument $v_i\in C^{\infty}(B_1)\cap C(\2{B})$ for all $i=1,2,\dots,m$. 
By (4), 
$$
v_i(x)=\Delta^{m-i}u(x)\quad\forall |x|=1,i=1,2,\dots,m.\tag 6
$$
By (3),
$$
\max_{|x|\le 1}|v_1|\le\max_{|x|=1}|\Delta^{m-1}u(x)|.\tag 7
$$
By (1) $\Delta^{m-1}u$ is harmonic in $B_1\setminus\{0\}$. Hence by (2), 
(5), (6), and standard theory of removable singularity for harmonic 
functions \cite{F},
$$
\Delta^{m-1}u(x)=v_1(x)\quad\text{ in }B_1\setminus\{0\}.\tag 8
$$
Hence $\Delta^{m-1}u$ has a removable singularity at $0$ and we can extend
$\Delta^{m-1}u$ to a smooth function on $B_1$ by letting $\Delta^{m-1}u(0)
=v_1(0)$. By (4) and (6), $\forall i=2,3,\dots,m$, 
$$\align
\|v_i\|_{L^{\infty}(B_1)}\le&C_1\|v_{i-1}\|_{L^{\infty}(B_1)}
+\|\Delta^{m-i}u\|_{L^{\infty}(\1 B_1)}\\
\le&C_1(C_1\|v_{i-2}\|_{L^{\infty}(B_1)}
+\|\Delta^{m-i+1}u\|_{L^{\infty}(\1 B_1)})
+\|\Delta^{m-i}u\|_{L^{\infty}(\1 B_1)}\\
\le&\cdots\\
\le&C_1^{i-1}\|\Delta^{m-1}u\|_{L^{\infty}(\1 B_1)}
+C_1^{i-2}\|\Delta^{m-2}u\|_{L^{\infty}(\1 B_1)}+\cdots
+\|\Delta^{m-i}u\|_{L^{\infty}(\1 B_1)}\tag 9
\endalign
$$
where
$$
C_1=\max_{|x|\le 1}\int_{B_1}G(x,y)\,dy=\max_{|x|\le 1}\frac{|x|^2}{2n}
=\frac{1}{2n}.
$$
We now claim that 
$$
\Delta^{m-i}u=v_i\quad\text{ in }B_1\setminus\{0\}\tag 10
$$
for any $i=1,2,\dots,m$.
We will prove the claim by induction. By the previous discussion the claim
holds for $i=1$. Suppose the claim holds for some $i=i_0\in\{1,2,\dots,
m-1\}$. Then by (4),
$$\align
&\Delta v_{i_0+1}=v_{i_0}=\Delta^{m-i_0}u\quad\text{ in }B_1\setminus\{0\}\\
\Rightarrow\quad&\Delta (\Delta^{m-i_0-1}u-v_{i_0+1})=0
\quad\text{ in }B_1\setminus\{0\}.\tag 11
\endalign
$$
By (2), (6), (9), (11) and standard theory of removable singularity for 
harmonic functions \cite{F} we get that (10) holds for $i=i_0+1$. Hence
by induction (10) holds for any $i=1,2,\dots,m$. Since $v_i\in C^{\infty}(B_1)$
for any $i=1,2,\dots,m$, by defining $\Delta^{m-i}u(0)=v_i(0)$ for all 
$i=1,2,\dots,m$, $\Delta^{m-i}u\equiv v_i\in C^{\infty}(B_1)$ for all 
$i=1,2,\dots,m$. Hence $u$ has a removable singularity at $x=0$.
\enddemo

\proclaim{\bf Theorem 2}
Let $x_0\in\Omega$. Suppose $u$ is a solution of (1) in $\Omega\setminus
\{x_0\}$ which satisfies
$$
|u(x)|=\left\{\aligned
&o(|x-x_0|^{2m-n})\qquad\qquad\qquad\quad\,\,\text{ as }|x-x_0|\to 0
\quad\text{ if }n\ge 3\\
&o(|x-x_0|^{2m-2}\log (|x-x_0|^{-1}))\quad\text{ as }|x-x_0|\to 0
\quad\text{ if }n=2.
\endaligned\right.
\tag 12
$$  
Then $u$ has a removable singularity at $x_0$.
\endproclaim
\demo{Proof}
With loss of generality we may assume that $x_0=0$ and $\2{B}_1\subset\Omega$. 
Since $u$ satisfies (1) in $\Omega\setminus\{0\}$, $\Delta^{m-1}u$ is 
harmonic in $\Omega\setminus\{0\}$. By the mean value theorem for harmonic
functions (cf. Appendix C of \cite{S}),
$$\align
\Delta^{m-1}u(x)=&(\phi_{\frac{|x|}{2}}\ast\Delta^{m-1}u)(x)
=((\Delta^{m-1}\phi_{\frac{|x|}{2}})\ast u)(x)\quad\forall 0<|x|\le 1/2\\
\Rightarrow\quad|\Delta^{m-1}u(x)|\le&\biggl (\frac{2}{|x|}\biggr )^n
\int_{|z|\le\frac{|x|}{2}}\biggl |\Delta_z^{m-1}\phi\biggl (
\frac{z}{|x|/2}\biggr )\biggr|
|u(x-z)|\,dz\\
\le&\biggl (\frac{2}{|x|}\biggr )^{n+2(m-1)}\int_{|z|\le\frac{|x|}{2}}
\biggl |(\Delta^{m-1}\phi)\biggl (\frac{z}{|x|/2}\biggr )\biggr ||u(x-z)|\,dz\\
\le&\frac{C}{|x|^{n+2(m-1)}}\int_{|z|\le\frac{|x|}{2}}
|u(x-z)|\,dz\quad\forall 0<|x|\le 1/2.\tag 13
\endalign
$$ 
By (12), for any $\3>0$ there exists $0<\delta_1<1$ such that
$$
|u(x)|\le\left\{\aligned
&\3|x|^{2m-n}\qquad\qquad\quad\forall 0<|x|\le\delta_1
\quad\text{ if }n\ge 3\\
&\3|x|^{2m-2}\log (1/|x|)\quad\forall 0<|x|\le\delta_1\quad\text{ if }n=2
\endaligned\right.\tag 14
$$
Then
$$\align
|u(x-z)|\le&\left\{\aligned
&\3|x-z|^{2m-n}\qquad\qquad\qquad\,\,\,\forall |z|\le\frac{|x|}{2},
0<|x|\le\frac{\delta_1}{2}\quad\text{ if }n\ge 3\\
&\3|x-z|^{2m-2}\log (1/|x-z|)\quad\forall |z|\le\frac{|x|}{2},
0<|x|\le\frac{\delta_1}{2}\quad\text{ if }n=2\endaligned\right.\\
\le&\left\{\aligned
&C_2\3|x|^{2m-n}\qquad\qquad\quad\forall |z|\le\frac{|x|}{2},
0<|x|\le\delta_2\quad\text{ if }n\ge 3\\
&C_2\3|x|^{2m-2}\log (1/|x|)\quad\forall |z|\le\frac{|x|}{2},
0<|x|\le\delta_2\quad\text{ if }n=2\endaligned\right.\tag 15
\endalign
$$
for some constants $0<\delta_2<\delta_1/2$ and $C_2>0$. By (13) and (15),
$$
|\Delta^{m-1}u(x)|\le\left\{\aligned 
&C_3\3|x|^{2-n}\qquad\,\,\,\,\forall 0<|x|\le\delta_2\quad\text{ if }n\ge 3\\
&C_3\3\log (1/|x|)\quad\forall 0<|x|\le\delta_2\quad\text{ if }n=2
\endaligned\right.
\tag 16
$$
for some constant $C_3>0$.
Let $v_1,v_2,\dots,v_m$ be given by (3) and (4). Then by (16) and an argument
similar to the proof of Theorem 1, (8) holds. Hence $\Delta^{m-1}u(x)$ has a 
removalbe singularity at $x=0$. By letting $\Delta^{m-1}u(0)=v_1(0)$, 
$\Delta^{m-1}u(x)$ is extended to a smooth function on $B_1$. 
We now claim that both (10) and
$$
|\Delta^{m-i}u(x)|=\left\{\aligned
&o(|x|^{2-n})\qquad\,\,\,\text{ as }|x|\to 0\quad\text{ if }n\ge 3\\
&o(\log (1/|x|))\quad\text{ as }|x|\to 0\quad\text{ if }n=2
\endaligned\right.
\tag 17
$$
holds for any $i=1,2,\dots,m$. We will prove the claim by induction. By the 
previous discussion the claim holds for $i=1$. Suppose (10) and (17) holds 
for some $i=i_0\in\{1,2,\dots,m-1\}$. Similar to the proof of Theorem 1 
(11) holds. Then by the mean value theorem for harmonic functions, 
$\forall 0<|x|\le 1/2$,
$$\align
\Delta^{m-i_0-1}u(x)-v_{i_0+1}(x)=&\phi_{\frac{|x|}{2}}\ast
(\Delta^{m-i_0-1}u-v_{i_0+1})(x)\\
=&\Delta^{m-i_0-1}\phi_{\frac{|x|}{2}}\ast u(x)-\phi_{\frac{|x|}{2}}\ast
v_{i_0+1}(x)\\
\Rightarrow\qquad\Delta^{m-i_0-1}u(x)=&\Delta^{m-i_0-1}\phi_{\frac{|x|}{2}}
\ast u(x)+v_{i_0+1}(x)-\phi_{\frac{|x|}{2}}\ast v_{i_0+1}(x).
\tag 18
\endalign
$$
Now by (15),
$$\align
|\Delta^{m-i_0-1}\phi_{\frac{|x|}{2}}\ast u(x)|
\le&\biggl (\frac{2}{|x|}\biggr )^n
\int_{|z|\le\frac{|x|}{2}}\biggl |\Delta_z^{m-i_0-1}\phi\biggl (
\frac{z}{|x|/2}\biggr )\biggr|
|u(x-z)|\,dz\\
\le&\biggl (\frac{2}{|x|}\biggr )^{n+2(m-i_0-1)}\int_{|z|\le\frac{|x|}{2}}
\biggl |(\Delta^{m-i_0-1}\phi)\biggl (\frac{z}{|x|/2}\biggr )
\biggr ||u(x-z)|\,dz\\
\le&\frac{C}{|x|^{n+2(m-i_0-1)}}\int_{|z|\le\frac{|x|}{2}}
|u(x-z)|\,dz\\
\le&\left\{\aligned
&C\3|x|^{2(i_0+1)-n}\qquad\,\,\forall 0<|x|\le\delta_2\text{ if }n\ge 3\\
&C\3|x|^{2i_0}\log (1/|x|)\quad\forall 0<|x|\le\delta_2\text{ if }n=2
\endaligned\right.\\
\le&\left\{\aligned
&C'\3|x|^{2-n}\qquad\,\,\,\forall 0<|x|\le\delta_2\text{ if }n\ge 3\\
&C'\3\log (1/|x|)\quad\forall 0<|x|\le\delta_2\text{ if }n=2.
\endaligned\right.\tag 19
\endalign
$$ 
By (9), (18) and (19), there exists a constants $C>0$ such that
$$
|\Delta^{m-i_0-1}u(x)|\le\left\{\aligned
&C\3|x|^{2-n}\qquad\,\,\,\forall 0<|x|\le\delta_2\text{ if }n\ge 3\\
&C\3\log (1/|x|)\quad\forall 0<|x|\le\delta_2\text{ if }n=2.
\endaligned\right.\tag 20
$$
Hence (17) holds 
for $i=i_0+1$. By (9), (11), (20) and standard theory for removable singularity
of harmonic functions, (10) holds for $i=i_0+1$. By induction (10) and (17)
holds for all $i=1,2,\dots,m$. Since $v_i\in C^{\infty}(B_1)$ for any 
$i=1,2,\dots,m$, by defining $\Delta^{m-i}u(0)=v_i(0)$ for all $i=1,2,\dots,
m$, $\Delta^{m-i}u\equiv v_i\in C^{\infty}(B_1)$ for all $i=1,2,\dots,m$.
Hence $u$ has a removable singularity at $x=0$.
\enddemo

\proclaim{\bf Remark 1}
Suppose $n\ge 3$ and either 
$$
n\text{ is odd and }m\ge 1
$$ 
or
$$
n\text{ is even and }m=1,\dots,\frac{n}{2}-1.
$$
By direct computation the function (cf. \cite{ACL}) $\Gamma (x)=|x|^{2m-n}$ 
satisfies (1) in $\Bbb{R}^n\setminus\{0\}$ and
$$
\Delta^{m-1}\Gamma(x)=A_{m,n}|x|^{2-n}\quad\forall |x|>0\tag 21
$$ 
for some constant $A_{m,n}\ne 0$. If $\Gamma(x)$ has a removable singularity 
at $x=0$, then $\Gamma(x)$ can be extended to a smooth solution of (1) in 
a small neighborhood of $0$. Hence $\Delta^{m-1}\Gamma (x)\in C^{\infty}$ 
in a small neighborhood of $0$. This contradicts (21). Hence $\Gamma(x)$ 
has a non-removable singularity at $x=0$. Thus Theorem 2 is sharp.
\endproclaim


\Refs

\ref
\key A\by D.H.~Armitage\paper On polyharmonic functions in $R^n
\setminus\{0\}$\jour J London Math. Soc. (2)\vol 8\yr 1974\pages 561--569
\endref

\ref
\key ACL\by \ \ \ N.~Aronszajn, T.M.~Creese and L.J.~Lipkin\book
Polyharmonic functions\publ Clarendon Press\publaddr Oxford\yr 1983
\endref

\ref
\key CGS\by \ \ \ L.A.~Cafferelli, B.~Gidas and J.~Spruck\paper
Asymptotic symmetry and local behaviour of semilinear elliptic 
equations with critical Sobolev growth\jour Comm. Pure Appl. Math.
\vol 42\yr 1989\pages 271--297\endref

\ref
\key F\by G.B.~Folland\book Introduction to Partial Differential 
Equations\publ Princeton University Press and University of
Tokyo Press\publaddr Princeton, New Jersey\yr 1976\endref

\ref 
\key H\by R.S.~Hamilton\paper The formation of singularities in the Ricci 
flow\jour Surveys in differential geometry, Vol. II (Cambridge, MA, 1993),
7--136, International Press, Cambridge, MA, 1995\endref

\ref
\key Hs1\by \ \ \ \ S.Y.~Hsu\paper Removable singularities and 
non-uniqueness of solutions of a singular diffusion equation 
\jour Math. Annalen\vol 325(4)\yr 2003\pages 665--693\endref

\ref
\key Hs2\by \ \ \ \ S.Y.~Hsu\paper A harmonic map flow associated with 
the standard solution of Ricci flow, http://arXiv.org\linebreak
/abs/math.DG/0702168\endref

\ref
\key P1\by G.~Perelman\paper The entropy formula for the Ricci flow and its 
geometric applications,  http://arXiv.org\linebreak /abs/math.DG/0211159
\endref 

\ref
\key P2\by G.~Perelman\paper Ricci flow with surgery on three-manifolds,
http://arXiv.org/abs/math.DG/0303109\endref

\ref
\key PV\by \ J.~Pipher and G.C.~Verchota\paper Maximum principles for the 
polyharmonic equation on Lipschitz domains\jour Potential Analysis\vol 4
\yr 1995\pages 615--636\endref

\ref
\key S\by E.M.~Stein\book Singular integral and differentiability properties
of functions\publ Princeton University Press\publaddr Princeton, New Jersey
\yr 1970\endref

\ref
\key X\by X.~Xu\paper Uniqueness theorem for the entire positive solutions
of biharmonic equations in $\Bbb{R}^n$\jour Proc. Roy. Soc. Edinburgh 
\vol 130A\pages 651--670\yr 2000\endref

\endRefs
\enddocument